\documentclass{amsart}
\usepackage{ae}
\usepackage{aecompl}
\usepackage[T1]{fontenc}
\usepackage[latin1]{inputenc}
\usepackage{amssymb}

 \theoremstyle{plain}
\newtheorem{thm}{Theorem}[section]
  \theoremstyle{definition}
  \newtheorem{defn}[thm]{Definition}
  \theoremstyle{plain}
  \newtheorem{lem}[thm]{Lemma}

\usepackage{mathbbol}

\usepackage[all,dvips,arc,curve,rotate]{xy}

\begin{document}

\title{Projective $C^{*}$-Algebras and Boundary Maps }

\author{Terry A. Loring}

\address{Department of Mathematics and Statistics, University of New Mexico,
Albuquerque, NM 87131, USA.}



\urladdr{http://www.math.unm.edu/\textasciitilde{}loring/}

\begin{abstract}
Both boundary maps in $K$-theory are expressed in terms of surjections
from projective $C^{*}$-algebras to semiprojective $C^{*}$-algebras. 
\end{abstract}

\maketitle

\markright {Projective $C^{*}$-Algebras and Boundary Maps (preprint version)}
\markleft {Projective $C^{*}$-Algebras and Boundary Maps (preprint version)}

\section{Noncommutative Cells and Boundaries}

Cells are absolute retracts that tie together spheres of different
dimensions. The analog of an absolute retract for a $C^{*}$-algebra
is being projective. For better or worse, in the category of all $C^{*}$-algebras,
we lose the projectivity of $C_{0}\left(\mathbb{D}\setminus\{-1\}\right),$
so we cannot generally use the exactness of
\[
0 \rightarrow C_{0}\left(\mathbb{R}^{2}\right) 
\rightarrow C_{0}\left(\mathbb{D}\setminus\{-1\}\right)
\rightarrow C_{0}\left(\mathbb{R}\right)\rightarrow0
\]
to explain the index map in $K$-theory. Another difficulty is that
we need asymptotic morphisms to obtain the natural isomorphism
\[
\left[\left[C_{0}\left(\mathbb{R}^{2}\right),D\otimes\mathbb{K}\right]\right]
\cong K_{0}(D).
\]

The name ``index map'' is related to the Toeplitz algebra $\mathcal{T}$
and the exact sequence
\[
0\rightarrow\mathbb{K}\rightarrow\mathcal{T}\rightarrow C(S^{1})\rightarrow0.
\]
We might prefer to use $\mathcal{T}_{0},$ generated by the shift
minus one, and
\[
0\rightarrow\mathbb{K}\rightarrow\mathcal{T}_{0}\rightarrow
C_{0}\left(\mathbb{R}\right)\rightarrow0,
\]
but still we may have trouble since $\mathcal{T}_{0}$ is not projective
and $\mathbb{K}$ is not semiprojective. 

The ``second standard picture of the index map'' in
\cite[Proposition 9.2.2]{RordamLandLbook}
and the picture of the exponential map presented in \cite{LoringProjectiveKtheory}
both use what might be called noncommutative cells. In both cases,
there is a diagram
\begin{equation}
\label{cellDiagram}
\xymatrix{
	0 \ar[r] &
		 U \ar[r] ^{\iota} &
			 P  \ar[r] ^{\eta} &
				 R  \ar[r]  &
					0 \\
	  & Q \ar[u] ^{\psi_0}
}
\end{equation}
with an exact row and where
$P$ is projective. Moreover, $R,$ $Q$ and $\psi_{0}$ have enough
nice properties to ensure that
\begin{equation}
\left[R,D\otimes\mathbb{K}\right]\cong K_{i}(D),\label{eq:K_i-isomorphism}
\end{equation}
\begin{equation}
\left[Q,D\otimes\mathbb{K}\right]\cong K_{i+1}(D),\label{eq:K_i+1-isomorphism}
\end{equation}
and $K_{i+1}(\psi_{0})$ is an isomorphism. The projectivity of $P$
then leads to an implementation of the boundary map as a sequence
of maps
\[
\partial^{(n)}:[R,\mathbf{M}_{n}(A/I)]\rightarrow[Q,\mathbf{M}_{n}(I)].
\]

There are other examples where we don't have the isomorphisms
(\ref{eq:K_i-isomorphism})
and (\ref{eq:K_i+1-isomorphism}). What we minimally require is the
following.

\begin{defn}
If the row in the diagram (\ref{cellDiagram}) is exact, $P$ is projective,
\[
K_{i}(R)=K_{i+1}(Q)=\mathbb{Z},\]
 \[
K_{i+1}(R)=K_{i}(Q)=0,\]
and $K_{i+1}(\psi_{0})$ is an isomorphism, then we will call (\ref{cellDiagram})
a \emph{cell diagram}. 
\end{defn}

\section{The Index Map\label{sec:The-Index-Map}}

The noncommutative Grassmannians are unital $C^{*}$-algebras with
universal properties. One gets easier statements of results if one
works with a nonunital variation. For now, we stick with the two-by-two
version, as this is the one most closely related to $q\mathbb{C}.$

We use $\tilde{A}$ to denote the unitization of $A,$ where a unit
$\mathbb{1}$ is always added.

For a set of relations $\mathcal{R}$ on a set $\mathcal{G}$ we use
the notation
\[
\iota:\mathcal{G}\rightarrow
C^{*}\left\langle \mathcal{G}\left|\mathcal{R}\right.\right\rangle
\]
to denote the function into a $C^{*}$-algebra that is the universal
representation of $\mathcal{R}.$ See \cite{Loring-lifting-perturbing}
for information on what relations are allowed. The definition of universal
representation requires that $\varphi\mapsto\varphi\circ\iota$ determines
a natural bijection between 
\[
\hom\left(C^{*}\left\langle
\mathcal{G}\left|\mathcal{R}
\right.\right\rangle ,B\right)
\]
and the set
\[
\left\{ f:\mathcal{G}\rightarrow B
\left|f\mbox{ is a representation of }\mathcal{R}
\right.\right\} .
\]

We can similarly work with relations in unital $C^{*}$-algebras,
and denote the universal representation in a unital $C^{*}$-algebra
by
\[
\iota:\mathcal{G}\rightarrow C_{1}^{*}
\left\langle \mathcal{G}\left|\mathcal{R}\right.
\right\rangle .
\]

In all the examples considered, $\iota$ will be an inclusion and
so we will identify $\mathcal{G}$ with $\iota(\mathcal{G}).$

Define $G_{2}^{\textnormal{nc}},$ c.f. \cite{Brown-free_product},
as
\[
G_{2}^{\textnormal{nc}}=
C_{1}^{*}\left\langle a,b,c\left|P^{2}=P^{*}=P\mbox{ for }
	P=\left[\begin{array}{cc}
	a & c^{*}\\
	c & b\end{array}\right]
\right.\right\rangle
\]
and define
\[
G_{2}^{\textnormal{st}}=
C^{*}\left\langle h,k,x\left|P^{2}=P^{*}=P\mbox{ for }
	P=\left[\begin{array}{cc}
	\mathbb{1}-h & x^{*}\\
	x & k\end{array}
\right]\right.\right\rangle .\]
The ``$\mbox{st}$'' is to stand for ``standard,'' as in the standard
picture of $K_{0}.$ The fact that a projection has norm at most one
means these relations are bounded, so this universal $C^{*}$-algebra
does exist.

\begin{lem}
The unitization $\left(G_{2}^{\textnormal{st}}\right)^{\sim}$ is
isomorphic to $G_{2}^{\textnormal{nc}}$ via
\begin{eqnarray*}
\mathbb{1} & \mapsto & 1\\
h & \mapsto & 1-a\\
k & \mapsto & b\\
x & \mapsto & c.
\end{eqnarray*}

\end{lem}

\begin{proof}
In terms of $*$-polynomial relations, $\left(G_{2}^{\textnormal{st}}\right)^{\sim}$
has generators $\mathbb{1},$ $h,$ $k,$ $x$ where $\mathbb{1}$
acts as a unit and
\begin{align}
 & h=h^{*}\nonumber \\
 & k=k^{*}\nonumber \\
 & h^{2}+x^{*}x=h\label{eq:GrassmanStandardRelations}\\
 & -xh+kx=0\nonumber \\
 & k^{2}+xx^{*}=k.\nonumber 
\end{align}
Clearly
\begin{eqnarray*}
h=h^{*} & \Longleftrightarrow & (\mathbb{1}-h)=(\mathbb{1}-h)^{*}\\
h^{2}+x^{*}x=h & \Longleftrightarrow & (\mathbb{1}-h)^{2}+x^{*}x=(\mathbb{1}-h)\\
-xh+kx=0 & \Longleftrightarrow & x(\mathbb{1}-h)+kx=x
\end{eqnarray*}
and
the result now follows easily.
\end{proof}

\begin{lem}
The $C^{*}$-algebra $G_{2}^{\textnormal{st}}$ is semiprojective. 
\end{lem}

\begin{proof}
This
follows from Corollary 2.16 and Proposition 2.17 of 
\cite{Blackadar-shape-theory}.
\end{proof}

Consider the automorphism
\[
\eta:G_{2}^{\textnormal{st}}\rightarrow G_{2}^{\mbox{st}}
\]
defined by $\eta(h)=k,$ $\eta(k)=h$ and $\eta(x)=x^{*}.$

\begin{lem}
The $*$-homomorphism
\[
\mbox{id}\oplus\eta:
G_{2}^{\textnormal{st}}\rightarrow
\mathbf{M}_{2}\left(G_{2}^{\textnormal{st}}\right)
\]
is null-homotopic.
\end{lem}

\begin{proof}
In terms of the generators, $h,$ $k$ and $x$ are being sent to
\[
\left[\begin{array}{cc}
h & 0\\
0 & k\end{array}\right],\quad\left[\begin{array}{cc}
k & 0\\
0 & h\end{array}\right],\quad\left[\begin{array}{cc}
x & 0\\
0 & x^{*}\end{array}\right].
\]
The homotopy is found in two segments.

For $0\leq\alpha<1,$ let $\beta=\sqrt{1-\alpha^{2}}.$ Let
\[
H_{\alpha}=\left[\begin{array}{cc}
h & 0\\
0 & k\end{array}\right],\quad K_{\alpha}=\left[\begin{array}{cc}
k & 0\\
0 & h\end{array}\right],
\]
\[
X_{\alpha}=\left[\begin{array}{cc}
\alpha x & -\beta\sqrt{xx^{*}}\\
\beta\sqrt{x^{*}x} & \alpha x^{*}\end{array}\right].
\]
Clearly $H_{\alpha}$ and $K_{\alpha}$ are self-adjoint. The commutation
relation is easy, since
\begin{eqnarray*}
 &  & \left[\begin{array}{cc}
\alpha x & -\beta\sqrt{xx^{*}}\\
\beta\sqrt{x^{*}x} & \alpha x^{*}\end{array}\right]\left[\begin{array}{cc}
h & 0\\
0 & k\end{array}\right]\\
 & = & \left[\begin{array}{cc}
\alpha kx & -\beta k\sqrt{k-k^{2}}\\
\beta h\sqrt{h-h^{2}} & \alpha hx^{*}\end{array}\right]\\
 & = & \left[\begin{array}{cc}
k & 0\\
0 & h\end{array}\right]\left[\begin{array}{cc}
\alpha x & -\beta\sqrt{xx^{*}}\\
\beta\sqrt{x^{*}x} & \alpha x^{*}\end{array}\right].
\end{eqnarray*}
For the remaining relations, we have
\begin{eqnarray*}
X_{\alpha}^{*}X_{\alpha} & = & \left[\begin{array}{cc}
\alpha^{2}x^{*}x+\beta^{2}x^{*}x 
	& \alpha\beta\left(-x^{*}\sqrt{xx^{*}}+\sqrt{x^{*}x}x^{*}\right)\\
\alpha\beta\left(-\sqrt{xx^{*}}x+x\sqrt{x^{*}x}\right) 
	& \alpha^{2}xx^{*}+\beta^{2}xx^{*}\end{array}\right]\\
 & = & \left[\begin{array}{cc}
h & 0\\
0 & k\end{array}\right]-\left[\begin{array}{cc}
h & 0\\
0 & k\end{array}\right]^{2}
\end{eqnarray*}
and by symmetry,
\[
X_{\alpha}X_{\alpha}^{*}=\left[\begin{array}{cc}
k & 0\\
0 & h\end{array}\right]-\left[\begin{array}{cc}
k & 0\\
0 & h\end{array}\right]^{2}.
\]
For the second part of the path, for each $0\leq\gamma\leq1,$ the
generators are
\[
H_{\gamma}=\left[\begin{array}{cc}
\gamma h & 0\\
0 & \gamma k\end{array}\right],\quad K_{\gamma}=\left[\begin{array}{cc}
\gamma k & 0\\
0 & \gamma h\end{array}\right],
\]
\[
X_{\gamma}=\left[\begin{array}{cc}
0 & -\sqrt{\gamma k-(\gamma k)^{2}}\\
\sqrt{\gamma h-(\gamma h)^{2}} & 0\end{array}\right].\]
Again the self-adjoint conditions are clear, and then
\begin{eqnarray*}
 &  & \left[\begin{array}{cc}
0 & -\sqrt{\gamma k-(\gamma k)^{2}}\\
\sqrt{\gamma h-(\gamma h)^{2}} & 0\end{array}\right]\left[\begin{array}{cc}
\gamma h & 0\\
0 & \gamma k\end{array}\right]\\
 & = & \left[\begin{array}{cc}
0 & -\gamma k\sqrt{\gamma k-(\gamma k)^{2}}\\
\gamma h\sqrt{\gamma h-(\gamma h)^{2}} & 0\end{array}\right]\\
 & = & \left[\begin{array}{cc}
\gamma k & 0\\
0 & \gamma h\end{array}\right]\left[\begin{array}{cc}
0 & -\sqrt{\gamma k-(\gamma k)^{2}}\\
\sqrt{\gamma h-(\gamma h)^{2}} & 0\end{array}\right]
\end{eqnarray*}
and
\begin{eqnarray*}
X_{\gamma}^{*}X_{\gamma} & = & \left[\begin{array}{cc}
\gamma h-(\gamma h)^{2} & 0\\
0 & \gamma k-(\gamma k)^{2}\end{array}\right]\\
 & = & \left[\begin{array}{cc}
\gamma h & 0\\
0 & \gamma k\end{array}\right]-\left[\begin{array}{cc}
\gamma h & 0\\
0 & \gamma k\end{array}\right]^{2}
\end{eqnarray*}
and by symmetry,\[
X_{\gamma}X_{\gamma}^{*}=\left[\begin{array}{cc}
\gamma k & 0\\
0 & \gamma h\end{array}\right]-\left[\begin{array}{cc}
\gamma k & 0\\
0 & \gamma h\end{array}\right]^{2}.\]

\end{proof}

The next result should be compared to the well-known isomorphisms
\[
\lim_{\rightarrow}\left[C_{0}(0,1),\mathbf{M}_{n}(D)\right]\cong K_{1}(D)
\]
and
\[
\lim_{\rightarrow}\left[q\mathbb{C},\mathbf{M}_{n}(D)\right]\cong K_{0}(D).
\]

\begin{thm}
For a $C^{*}$-algebra $D,$ there is a natural isomorphism
\[
\lim_{\rightarrow}\left[G_{2}^{\textnormal{st}},\mathbf{M}_{n}(D)\right]
\cong K_{0}(D).
\]

\end{thm}

\begin{proof}
By \cite[Theorem 4.3]{DadarlatLoringUnsuspendedE} we know
\[
K_{0}(D)\cong\left[\left[G_{2}^{\mbox{st}},A\otimes\mathbb{K}\right]\right]
\]
and by semiprojectivity
\[
\left[\left[G_{2}^{\mbox{st}},D\otimes\mathbb{K}\right]\right]
\cong
\left[G_{2}^{\mbox{st}},D\otimes\mathbb{K}\right]
\cong
\lim_{\rightarrow}\left[G_{2}^{\mbox{st}},D\otimes\mathbf{M}_{n}\right].
\]

This also follows from standard results in $K$-theory.
\end{proof}

Recall from \cite{Cuntz-new-look-KK} that $q\mathbb{C}$ was defined
via an exact sequence
\[
0
\rightarrow q\mathbb{C}
\rightarrow\mathbb{C}\ast\mathbb{C}
\rightarrow\mathbb{C}
\rightarrow0.
\]
It has the concrete description
\[
q\mathbb{C}=
\left\{
f\in C_{0}\left((0,1],\mathbf{M}_{2}\right)\left|f(1)\mbox{ is diagonal}\right.
\right\} .
\]
as well as being universal on generators $h_{0},$ $k_{0}$ and $x_{0}$
for the relations
\[
\mathcal{P}=
C^{*}\left\langle h_{0},k_{0},x_{0}\left|h_{0}k_{0}=0,
\ P_{0}^{2}=P_{0}^{*}=P_{0}\mbox{ for }
	P_{0}=\left[\begin{array}{cc}
	\mathbb{1}-h_{0} & x_{0}^{*}\\
	x_{0} & k_{0}
	\end{array}\right]
\right.\right\rangle .\]

\begin{thm}
\label{thm:G_nc^st-KK-to-qC}
There is a surjection $\rho$ and an
inclusion $\lambda$ so that
\[
\xymatrix{
\mathbf {M}_{2}\left(G_{2}^{\textnormal{st}}\right ) \\
G_{2}^{\textnormal{st}} 
	\ar[u] ^{\textnormal{id} \otimes e_{11}}
	\ar[r] ^{ \rho}
 		& q\mathbb{C} \ar[lu] _{\lambda }
}
\mbox{ and }
\xymatrix{
\mathbf {M}_{2}\left(q\mathbb{C}\right ) \\
q\mathbb{C}
	\ar[u] ^{\textnormal{id} \otimes e_{11}}
	\ar[r] ^(0.45){ \lambda}
 		& \mathbf {M}_{2}\left( G_{2}^{\textnormal{st}} \right ) 
			\ar[lu] _{\rho \otimes  \textnormal{id}}
}
\]
commute up to homotopy.
In terms of generators,\[
\rho(h)=h_{0},\]
\[
\rho(k)=k_{0},\]
\[
\rho(x)=x_{0}\]
and\[
\lambda(h_{0})=h\otimes e_{11},\]
\[
\lambda(k_{0})=k\otimes e_{22},\]
\[
\lambda(x_{0})=x\otimes e_{21}.\]
Composition with $\lambda$ leads to a natural isomorphism\[
\lim_{\rightarrow}\left[G_{2}^{\mbox{st}},D\otimes\mathbf{M}_{n}\right]
\cong
\lim_{\rightarrow}\left[q\mathbb{C},D\otimes\mathbf{M}_{n}\right]
\]
(and with $K_{0}(D).$)
\end{thm}
\begin{proof}
We can define the homotopy $\varphi_{t}$ from
$\lambda\circ\rho$ to $\mbox{id}\otimes e_{11}$ on generators by
\[
\varphi_{t}(h)=h\otimes|w_{t}|,
\]
\[
\varphi_{t}(k)=k\otimes|w_{t}^{*}|,
\]
\[
\varphi_{t}(x)=x\otimes w_{t}\]
for some homotopy of partial isometries $w_{t}$ from $e_{21}$ to
$e_{11}.$ The homotopy from $(\rho\otimes\mbox{id})\circ\lambda$
to $ $$\mbox{id}\otimes e_{11}$ is found in a similar manner.
\end{proof}

Now we look at an extension that is somehow universal for the index
map. See \cite{RordamLandLbook} and ``the second standard picture
of the index map,'' (proposition 9.2.2). Recall
\[
C_{0}(0,1)=
C^{*}\left\langle
x\left|\strut(\mathbb{1}+x)^{*}=(\mathbb{1}+x)^{-1}\right.
\right\rangle 
\]
and define
\[
\mathcal{D}=
C^{*}\left\langle y\left|\strut\|\mathbb{1}+y\|\leq1\right.
\right\rangle .
\]
Sending $y$ to $x$ gives a surjection. Let $\mathcal{V}$ be the
kernel, so that we have the exact sequence 
\[
\xymatrix{
	0 \ar[r] &
		\mathcal V \ar[r] &
			\mathcal D  \ar[r] &
				C_0(0,1)  \ar[r] &
					0
} .
\]

\begin{lem}
The $C^{*}$-algebra $\mathcal{D}$ is projective and $C_{0}(0,1)$
is semiprojective.
\end{lem}
\begin{proof}
These have unitization the universal unital $C^{*}$-algebra generated
by a contraction and $C(S^{1}),$ respectively. The usual facts about
unitaries and contractions tell us these are semiprojective, or in
the first case projective, in the unital category. We are done, by
Theorem 10.1.9 and Lemma 14.1.6 of \cite{Loring-lifting-perturbing}.
\end{proof}
Consider $a=\mathbb{1}+x$ and\begin{eqnarray*}
h_{1} & = & \mathbb{1}-a^{*}a\\
k_{1} & = & \mathbb{1}-aa^{*}\\
x_{1} & = & a\sqrt{\mathbb{1}-a^{*}a}.\end{eqnarray*}
These are all elements of $\mathcal{V}$ and it is easy to see that
\[
P_{1}=\left[\begin{array}{cc}
\mathbb{1}-h_{1} & x_{1}^{*}\\
x_{1} & k_{1}\end{array}\right]
\]
is a projection. 
This determines a $*$-homomorphism $\psi_{0}$ from $G_{2}^{\mbox{st}}$
to $\mathcal{V}.$

\begin{lem}
The diagram 
\[
\xymatrix{
	0 \ar[r] &
		\mathcal V \ar[r] &
			\mathcal D  \ar[r] &
				 C_0(0,1)  \ar[r] &
					0 \\
	  & G_{2}^{\textnormal{st}} \ar[u] ^{\psi_0}
}
\]
is a cell diagram.
\end{lem}

\begin{proof}
Since $K_{*}(\mathcal{D})=0$ we know 
\[
\partial:K_{1}(C_{0}(0,1))\rightarrow K_{0}(\mathcal{V})
\]
is an isomorphism. The fact that $K_{1}(\psi_{0})$ is an isomorphism
follows from the definition of the boundary map (\cite{RordamLandLbook}).
\end{proof}

\section{The exponential Map\label{sec:The-exponential-Map}}

Consider a short exact sequence 
\[
0\rightarrow I\rightarrow A\rightarrow A/I\rightarrow0
\]
and the associated boundary map 
\[
\partial:K_{0}(A/I)\rightarrow K_{1}(I).
\]
In \cite{LoringProjectiveKtheory} we showed that if $x$ in $K_{0}(A/I)$
is realized by $\varphi$ in $\hom(q\mathbb{C},A/I)$ then $\partial(x)$
is realized by some $\psi$ in $\hom(C_{0}(0,1),I).$ Equivalently,
$\partial(x)$ is realized as a unitary in $\tilde{I}.$ 

We want to verify that the construction of $\psi$ is well defined
up to homotopy and that 
\[
[q\mathbb{C},A/I]\rightarrow[C_{0}(0,1),I]
\]
is natural. This can be done by an examination of the proof of 
\cite[Theorem 6]{LoringProjectiveKtheory},
but is more clearly taken care of by a cell diagram and 
Theorem~\ref{thm:projectiveCreatesBoundary}.

In this approach to the exponential map, the key point is the projectivity
of
\[
\mathcal{P}=
C^{*}\left\langle
h,k,x\left|hk=0,\ 0\leq P\leq1\mbox{ for }
P=\left[\begin{array}{cc}
\mathbb{1}-h & x^{*}\\
x & k\end{array}\right]
\right.
\right\rangle .
\]
We will not reprove that, but refer the reader to 
\cite[Theorem 9]{LoringProjectiveKtheory}.
What we will do is give a second approach to the $K$-theory calculations
related to $\mathcal{P}$ based on finding an embedding
\[
\mathcal{P}\hookrightarrow\mathbb{C}\ast C_{0}(0,1].
\]

Let $\eta$ be surjection $\eta:\mathcal{P}\rightarrow q\mathbb{C}$
onto\[
q\mathbb{C}=
C^{*}\left\langle
h_{0},k_{0},x_{0}\left|P_{0}^{*}=P_{0}^{2}=P_{0}\mbox{ for }P_{0}=\left[\begin{array}{cc}
\mathbb{1}-h_{0} & x_{0}^{*}\\
x_{0} & k_{0}\end{array}\right]\right.
\right\rangle 
\]
defined by $\eta(h)=h_{0},$ etc. Let $\mathcal{U}$ denote the kernel
of $\eta$ and $\iota$ denote the inclusion. Recall that $K_{0}(q\mathbb{C})$
is a copy of $\mathbb{Z}$ generated by the class of the projection
$P_{0}$ in $\mathbf{M}_{2}\left(\left(q\mathbb{C}\right)^{\sim}\right).$

\begin{lem}
\label{lem:ExponentialCell} There is a $*$-homomorphism $\psi_{0}$
so that
\[
\xymatrix{
0  \ar[r] & 
	\mathcal {U}  \ar@{^{(}->}[r] ^{\iota}&
		\mathcal {P}  \ar[r] ^{\eta}&
			q \mathbb {C}  \ar[r] &
			0 \\
			&
		C_0(0,1) \ar[u] ^{\psi_0}
}
\]
is a cell diagram.
\end{lem}

Since we know from \cite{LoringProjectiveKtheory} that $\mathcal{P}$
is projective, we need only find $\psi_{0}$ 
that induces an isomorphism on $K$-theory. This is easily seen to
be equivalent to the following Lemma.

\begin{lem}
\label{lem:exponentialFormula}In $\tilde{\mathcal{P}},$ let\[
u=-\mathbb{1}+\sum_{i,j}v_{ij}\]
where\[
v=\left[\begin{array}{cc}
v_{11} & v_{12}\\
v_{21} & v_{22}\end{array}\right]=e^{2\pi iP}.\]
Then $u$ is a unitary in $\tilde{U}$ that represents $\partial\left(\left[P_{0}\right]\right)$
in $K_{1}(I).$ 
\end{lem}

\begin{proof}
Theorem 6 in \cite{LoringProjectiveKtheory}, applied to
\[
0\rightarrow\mathcal{U}\rightarrow\mathcal{P}\rightarrow q\mathbb{C}\rightarrow0,
\]
tells us that some unitary in $\tilde{U}$ will represent that $K_{1}$-class
of the boundary $[P_{0}].$ The proof of that result
tells gives us the formula for $u,$ and moreover shows that $e^{2\pi iP_{0}}$
is homotopic through unitaries to 
\[
\left[\begin{array}{cc}
u & 0\\
0 & \mathbb{1}\end{array}\right].
\]
 
\end{proof}
The rest of this section is devoted to an alternative proof of 
Lemma~\ref{lem:ExponentialCell}.

First we get more specific regarding the exact sequence 
\[
\xymatrix{
0 \ar[r] & 
	q \mathbb {C}  \ar[r]^(0.45){\theta_0} &
		\mathbb {C} \ast \mathbb {C}  \ar[r] ^(0.6){\rho_0} &
			\mathbb {C}  \ar[r] &
				0 
} .
\]
We will use $p_{0}$ and $q_{0}$ to denote the two generating projections
in $\mathbb{C}\ast\mathbb{C}.$ Both of these are sent to $1$ by
$\rho_{0}.$ The inclusion $\theta_{0}$ of $q\mathbb{C}$ in $\mathbb{C}\ast\mathbb{C}$
is determined on generators by
\begin{eqnarray*}
\theta(h_{0}) & = & p_{0}-p_{0}q_{0}p_{0},\\
\theta(k_{0}) & = & (\mathbb{1}-p_{0})q_{0}(\mathbb{1}-p_{0}),\\
\theta(x_{0}) & = & (\mathbb{1}-p_{0})q_{0}p_{0}.
\end{eqnarray*}

There is a similar exact sequence involving $\mathcal{P}.$ Let the
obvious generators of $\mathbb{C}\ast C_{0}(0,1]$ be denoted $p$
and $l,$ so the only relations on them are
\begin{align*}
 & p^{2}=p^{*}=p,\\
 & 0\leq l\leq1.
 \end{align*}

\begin{thm}
\label{thm:PasIdealInFreeProduct}There are $*$-homomorphisms $\theta$
and $\rho$ defined by 
\begin{eqnarray*}
\theta(h) & = & p-plp\\
\theta(k) & = & (\mathbb{1}-p)l(\mathbb{1}-p)\\
\theta(x) & = & (\mathbb{1}-p)lp
\end{eqnarray*}
and
\begin{eqnarray*}
\rho(p) & = & 1,\\
\rho(l) & = & 1
\end{eqnarray*}
so that the sequence 
\[
\xymatrix{
0 \ar[r] & 
	\mathcal {P} \ar[r] ^(0.3){\theta} &
		\mathbb {C} \ast C_0 (0,1] \ar[r] ^(0.7){\rho} &
			\mathbb {C} \ar[r] &
				0
}
\]
is exact.
\end{thm}
\begin{proof}
Since
\begin{eqnarray*}
 &  & \left[\begin{array}{cc}
\mathbb{1}-\left(p-plp\right) & \left((\mathbb{1}-p)lp\right)^{*}\\
(\mathbb{1}-p)lp & (\mathbb{1}-p)l(\mathbb{1}-p)\end{array}\right]\\
 & = & \left[\begin{array}{cc}
p & \mathbb{1}-p\\
\mathbb{1}-p & p\end{array}\right]\left[\begin{array}{cc}
l & 0\\
0 & \mathbb{1}-p\end{array}\right]\left[\begin{array}{cc}
p & \mathbb{1}-p\\
\mathbb{1}-p & p\end{array}\right]
\end{eqnarray*}
we have
\[
0\leq\left[\begin{array}{cc}
\mathbb{1}-\left(p-plp\right) & \left((\mathbb{1}-p)lp\right)^{*}\\
(\mathbb{1}-p)lp & (\mathbb{1}-p)l(\mathbb{1}-p)\end{array}\right]\leq1.
\]
Since

\[
\left(p-plp\right)\left((\mathbb{1}-p)l(\mathbb{1}-p)\right)=0,
\]
we see that $\theta$ is well-defined. 

The unit $1$ is both a projection and a positive contraction, so
$\rho$ is well-defined.

Exactness at $\mathbb{C}$ is obvious.

To prove exactness at $\mathcal{P},$ suppose 
$\pi:\mathcal{P}\rightarrow\mathbb{B}(\mathbb{H})$
is a faithful representation of $\mathcal{P},$ and let $h_{1}=\pi(h),$
etc. Let $r=[h_{1}]$ be the range projection of $h_{1}$ and let
$q=[k_{1}]$ be the range projection of $k_{1}.$ The orthogonality
of $h$ and $k$ implies orthogonality for $r$ and $q.$ We established
in the proof of Theorem 4.3 of \cite{LoringProjectiveKtheory} the
factorization $x=k^{\frac{1}{8}}yh^{\frac{1}{8}}$ for some $y$ and
so
\[
rx_{1}=x_{1}q=0
\]
and
\[
qx_{1}=x_{1}r=x_{1},
\]
and of course
\[
rh_{1}=h_{1}r=h_{1},
\]
\[
qh_{1}=h_{1}q=rk_{1}=k_{1}r=0,
\]
\[
qk_{1}=k_{1}q=k_{1}.
\]
We know
\[
0\leq\left[\begin{array}{cc}
I-h_{1} & x_{1}^{*}\\
x_{1} & k_{1}\end{array}\right]\leq1
\]
and so
\[
0\leq\left[\begin{array}{cc}
r & q\end{array}\right]\left[\begin{array}{cc}
I-h_{1} & x_{1}^{*}\\
x_{1} & k_{1}\end{array}\right]\left[\begin{array}{cc}
r & q\end{array}\right]^{*}\leq1.
\]
This says
\[
0\leq r-h_{1}+x_{1}+x_{1}^{*}+k_{1}\leq1
\]
We can define a representation $\overline{\pi}$ of $\mathbb{C}\ast C_{0}(0,1]$
on $\mathbb{B}(\mathbb{H})$ by setting $\overline{\pi}(p)=r$ and
\[
\overline{\pi}(l)=r-h_{1}+x_{1}+x_{1}^{*}+k_{1}.
\]
This is an extension of $\pi$ because
\[
\overline{\pi}\circ\theta(h)=r-r(r-h_{1}+x_{1}+x_{1}^{*}+k_{1})r=h_{1}
\]
and
\[
\overline{\pi}\circ\theta(k)=(I-r)(r-h_{1}+x_{1}+x_{1}^{*}+k_{1})(I-r)=k_{1}
\]
and
\[
\overline{\pi}\circ\theta(x)=(I-r)(r-h_{1}+x_{1}+x_{1}^{*}+k_{1})r=x_{1}.
\]
We have shown $\iota$ is one-to-one, and so have exactness at $\mathcal{P}.$

Next we show that the image of $\theta$ is an ideal. This follows
from these equalities:
\begin{align*}
\left(p-plp\right)p 
	& =  \left(p-plp\right)\\
\left(p-plp\right)l 
	& =  \left(p-plp\right)\left((\mathbb{1}-p)lp\right)^{*}+\left(p-plp\right)
		-\left(p-plp\right)^{2}\\
\left((\mathbb{1}-p)l(\mathbb{1}-p)\right)p 
	& =  0\\
\left((\mathbb{1}-p)l(\mathbb{1}-p)\right)l 
	& =  \left((\mathbb{1}-p)l(\mathbb{1}-p)\right)\left((\mathbb{1}-p)lp\right)
		+\left((\mathbb{1}-p)l(\mathbb{1}-p)\right)^{2}\\
\left((\mathbb{1}-p)lp\right)p 
	& =  \left((\mathbb{1}-p)lp\right)\\
p\left((\mathbb{1}-p)lp\right)
	& =  0\\
\left((\mathbb{1}-p)lp\right)l 
	& =  \left((\mathbb{1}-p)lp\right)-
		\left((\mathbb{1}-p)lp\right)\left(p-plp\right)
			+\left((\mathbb{1}-p)lp\right)\left((\mathbb{1}-p)lp\right)^{*}\\
l\left((\mathbb{1}-p)lp\right) 
	& =  \left((\mathbb{1}-p)lp\right)^{*}\left((\mathbb{1}-p)lp\right)
	+\left((\mathbb{1}-p)l(\mathbb{1}-p)\right)\left((\mathbb{1}-p)lp\right).
\end{align*}

As to exactness in the middle, it is clear that $\rho\circ\theta=0.$
We need to show that the induced map 
\[
\overline{\rho}:
\left.\left(\mathbb{C}\ast C_{0}(0,1]\right)\right/\theta(\mathcal{P})
\rightarrow
\mathbb{C}
\]
is an isomorphism. Since
\begin{eqnarray*}
l-p 
	& = & -\left(p-plp\right)
	+\left((\mathbb{1}-p)lp\right)+\left(pl(\mathbb{1}-p)\right)
	+\left((\mathbb{1}-p)l(\mathbb{1}-p)\right)\\
	& = & -\theta(h)+\theta(x)+\theta(x)^{*}+\theta(k)
\end{eqnarray*}
we discover
\[
\left.\left(\mathbb{C}\ast C_{0}(0,1]\right)\right/\theta(\mathcal{P})
\]
is generated by a single projection. Since $\overline{\rho}$ maps
onto $\mathbb{C},$ it must be an isomorphism.
\end{proof}
Recall we have the surjection $\eta:\mathcal{P}\rightarrow q\mathbb{C}.$
Consider also the surjection 
\[
\eta_{1}:\mathbb{C}\ast C_{0}(0,1]\rightarrow\mathbb{C}\ast\mathbb{C}
\]
defined via
\begin{eqnarray*}
\eta_{1}(p) & = & p,\\
\eta_{1}(l) & = & q.
\end{eqnarray*}
Let us use $\mathcal{U}_{1}$ to denote the kernel of $\eta_{1}.$
This gives us the diagram
\[
\xymatrix{
	&	& 0 
			& 0
\\
	&	& \mathbb{C} \ar[u] \ar@{=}[r]
			& \mathbb{C} \ar[u]
\\
0 \ar[r] & 
	\mathcal {U}_1 \ar@{^{(}->}[r] ^(0.35){\iota_1}&
		\mathbb {C} \ast C_0 (0,1] \ar[r] ^(0.6){\eta_1} \ar[u] ^{\rho}&
			\mathbb {C} \ast \mathbb {C} \ar[r] \ar[u] ^{\rho_0} &
			0 \\
0  \ar[r] & 
	\mathcal {U} \ar[u] _{\cong}^{\theta_1} \ar@{^{(}->}[r] ^{\iota} &
		\mathcal {P} \ar@{^{(}->}[u] ^{\theta} \ar[r] ^{\eta}&
			q \mathbb {C} \ar@{^{(}->}[u] ^{\theta_0} \ar[r] &
			0 \\
	&	& 0 \ar[u]
			& 0		\ar[u]	
}
\]
where $\theta_{1}$
is the restriction of $\theta.$ Both rows and both columns are exact,
and it follows that $\theta_{1}$ is an isomorphism. 

The $K$-theory of the middle row is easy to work out, and so we see
that $K_{1}(\mathcal{U}_{1})\cong\mathbb{Z}$ and has generator represented
by the unitary $e^{2\pi il}$ in $\mathcal{U}_{1}^{\sim}.$ This completes
the second proof of Lemma~\ref{lem:ExponentialCell}.

\section{Projectives Determine Boundary Maps}

\begin{lem}
\label{lem:Lifting homotopies}Suppose $P$ is projective, $\rho:A\rightarrow B$
is a surjective $*$-homomorphism and $\varphi_{t}:P\rightarrow B$
is a homotopy of $*$-homomorphisms. Given $*$-homomorphisms $\psi_{0}$
and $\psi_{1}$ from $P$ to $A$ that are lifts of $\varphi_{0}$
and $\varphi_{1},$ there exists a homotopy of $*$-homomorphisms
$\bar{\psi}$ so that $\bar{\psi}_{t}$ is a lift of $\varphi_{t}$
and $\bar{\psi}_{0}=\psi_{0}$ and $\bar{\psi}_{1}=\psi_{1}.$ 
\end{lem}
\begin{proof}
This proof is a standard argument using a mapping cylinder.
\end{proof}

\begin{thm}
\label{thm:projectiveCreatesBoundary} Suppose
\[
\xymatrix{
	0 \ar[r] &
		 U \ar[r] ^{\iota} &
			 P  \ar[r] ^{\eta} &
				 R  \ar[r]  &
					0 \\
	  & Q \ar[u] ^{\psi_0}
}
\]
is a cell diagram. Suppose that $\alpha$ is a fixed generator of
$K_{i}(R)$ and that $\beta$ in $K_{i+1}(Q)$ defined so that 
$\partial(\alpha)=(\psi_{0})_{*}(\beta).$
Given an ideal $I$ in $A,$ there are natural maps $\partial^{(n)}$
so that 
\[
\xymatrix@R=1.3cm@C=1.3cm{
\left [  R, \mathbf {M}_{n} ( A/I ) \right ]  
				\ar[r] ^(0.54){\partial^{(n)}}
				\ar[d] \ar@/_0.45cm/[dd] |*[*3]{\hole}  _(0.25){\Xi^{(n)}}
	& \left [ Q, \mathbf {M}_{n} ( I ) \right ] 
				\ar[d] \ar@/_0.45cm/[dd] |*[*3]{\hole}  _(0.25){\Lambda^{(n)}}\\
\left [  R,\mathbf {M}_{n+1} ( A/I ) \right ]  
				\ar[r] ^(0.54){\partial^{({n+1})}}
				\ar[d]^(0.45){\Xi^{(n+1)}}
	& \left [ Q, \mathbf {M}_{n+1} ( I )  \right ] 
				\ar[d] ^(0.45){\Lambda^{(n+1)}}\\
K_i(A/I) \ar[r] ^{\partial}
	&  K_{i+1} (I)
}
\]
commutes, where $\partial$ is the boundary
map in $K$-theory. Here $\Xi^{(n)}(\psi)=\psi_{*}(\alpha)$
and $\Lambda^{(n)}(\varphi)=\varphi_{*}(\beta).$ 
\end{thm}

\begin{proof}
The naturality of $\partial^{(1)}$ allows us construct $\partial^{(n)}$
out of $\partial^{(1)}$ so we only concern ourselves with finding
$\partial^{(1)}$ with $\delta\circ\Xi^{(1)}=\Lambda^{(1)}\circ\delta^{(1)}.$

Suppose 
\[
\xymatrix{
0  \ar[r] & 
	I \ar[r] ^{\kappa}&
		A \ar[r] ^{\pi} &
			B   \ar[r] &
				0 
}
\]
is exact and we are given $\varphi:R\rightarrow B.$
Projectivity tells us there exists $\bar{\varphi}$ 
with $\pi\circ\bar{\varphi}=\varphi\circ\theta$
which we call $\bar{\varphi}.$ This restricts to a map $\hat{\varphi}$
between ideals. As a commuting diagram with exact rows, we are here: 
\[
\xymatrix{
	0 \ar[r]&
		I \ar[r] ^{\kappa} &
			A \ar[r] ^{\pi} &	
				B \ar[r] &
					0\\
0  \ar[r] & 
	U \ar[r] ^{\iota}  \ar[u] ^{\hat\varphi}&
		P  \ar[r] ^{\theta} \ar[u] ^{\bar\varphi}&
			R   \ar[r] \ar[u] _{\varphi}&
			0 \\
	& Q \ar[u] _{\psi_0}
}
\]

We wish to define
\[
\partial^{(1)}([\varphi])=\left[\hat{\varphi}\circ\psi_{0}\right].
\]
This composition depends on our choice of $\bar{\varphi}$ as well
as the choice of representative of the homotopy class $[\varphi].$
By the Lemma \ref{lem:Lifting homotopies}, we get
a well defined map from $[R,B]$ to $[Q,I].$ The naturality of the
boundary map in $K$-theory shows $\partial^{(1)}$ implements the
boundary map. The flexibility in choosing the lift $\bar{\varphi}$
makes it easy to show that $\partial^{(1)}$ is natural.
\end{proof}

Applying this to the examples in Sections~\ref{sec:The-Index-Map}
and \ref{sec:The-exponential-Map} we get the  somewhat unified picture
of the boundary maps summarized in three diagrams: 
\[
\xymatrix{
\left [  C_0(0,1),\mathbf {M}_{n} ( A/I ) \right ]  
				\ar[r] ^(0.54){\partial^{({n})}}
				\ar[d]
	& \left [ G_{2}^{\textnormal{st}}, \mathbf {M}_{n} ( I )  \right ] 
				\ar[d] 
				\\
\displaystyle{\lim_{\rightarrow}
\left [  C_0(0,1),\mathbf {M}_{k} ( A/I ) \right ]  
				\ar[r] 
				\ar[d]^{\cong}
}
				& \displaystyle{\lim_{\rightarrow}
				   \left [ G_{2}^{\textnormal{st}}, \mathbf {M}_{k} ( I )  \right ] 
					\ar[d] ^{\cong}
					}
				\\
K_1(A/I) \ar[r] ^{\partial}
	&  K_0 (I)
}
\]
\[
\xymatrix{
\left [  G_{2}^{\textnormal{st}},\mathbf {M}_{n} ( D ) \right ] 
				\ar[d] \ar[r] 
	& \left [ q \mathbb{C}, \mathbf {M}_{2n} ( D )  \right ] 
				\ar[d]
\\
\displaystyle{\lim_{\rightarrow}
\left [  G_{2}^{\textnormal{st}},\mathbf {M}_{k} ( D ) \right ] 
				\ar[d] ^{\cong}
				\ar[r] ^{\cong}
}
	& \displaystyle{\lim_{\rightarrow}
	\left [ q \mathbb{C}, \mathbf {M}_{k} ( D )  \right ] 
				\ar[dl]  ^{\cong}
	}
	\\
  K_0 (D)
}
\]
\[
\xymatrix{
\left [  q \mathbb{C},\mathbf {M}_{n} ( A/I ) \right ]  
				\ar[r] ^(0.54){\partial^{({n})}}
				\ar[d]
	& \left [ C_0(0,1), \mathbf {M}_{n} ( I )  \right ] 
				\ar[d] ^{\Lambda^{(n)}}\\
\displaystyle{\lim_{\rightarrow}
\left [  q \mathbb{C},\mathbf {M}_{k} ( A/I ) \right ]  
				\ar[d] ^{\cong}
				\ar[r]
}
	& \displaystyle{\lim_{\rightarrow}
	\left [ C_0(0,1), \mathbf {M}_{k} ( I )  \right ] 
				\ar[d] ^{\cong}
				}
				\\
K_0(A/I) \ar[r] ^{\partial}
	&  K_1 (I)
}
\]
The jump up in the size
of matrices in the middle diagram is the one we saw in
Theorem~\ref{thm:G_nc^st-KK-to-qC}.

\section{Further Examples}

The list of projective $C^{*}$-algebras is still growing.
For example, Shulman (\cite{Shulman_nilpotents}) has recently
shown that a nilpotent
contraction lifts to a nilpotent contraction of the same order.
It should be fruitful to search for semiprojective
quotients of projective $C^{*}$-algebras
that can also serve as ``boundaries of nonstandard cells.''
To illustrate, we now look at two more applications of 
Theorem~\ref{thm:projectiveCreatesBoundary}.

Recall that the cone\[
\mathbf{C}\mathbf{M}_{n}=C_{0}\left((0,1],\mathbf{M}_{n}\right)\]
is projective and fits in a nice exact sequence
\[
0\rightarrow\mathbf{S}\mathbf{M}_{n}
\rightarrow\mathbf{C}\mathbf{M}_{n}\rightarrow\mathbf{M}_{n}\rightarrow0
\]
with the suspension $\mathbf{S}\mathbf{M}_{n}$ of $\mathbf{M}_{n}.$
The larger cone\[
\mathbf{C}\mathbb{K}=C_{0}\left((0,1],\mathbb{K}\right)\]
is not residually finite dimensional, so it is not projective. We
can instead consider the mapping telescope
\[
\mathbf{T}(\mathbb{K})
=
\left\{
f\in C_{0}\left((0,\infty],\mathbb{K}\right)
\left|
\ t\leq n\implies f(t)\in\mathbf{M}_{n}
\right.
\right\} 
\]
which is projective (\cite{LorPederProjTrans}) and maps onto $\mathbb{K}$
by evaluation at $\infty.$ The kernel of this map is a bit awkward,
being
\[
\mathbf{I}(\mathbb{K})
=
\left\{
f\in C_{0}\left((0,\infty),\mathbb{K}\right)
\left|
\ t\leq n\implies f(t)\in\mathbf{M}_{n}
\right.
\right\} .
\]

\begin{thm}
Let $x$ be the standard generator of $K_{0}\left(\mathbf{M}_{n}\right)$
and $y$ be the standard generator of $K_{1}\left(\mathbf{S}\mathbf{M}_{n}\right).$
Given an ideal $I$ in a $C^{*}$-algebra $A,$ for any $*$-homomorphism
\[
\varphi:\mathbf{M}_{n}\rightarrow A/I
\]
there is a $*$-homomorphism
\[
\psi:\mathbf{S}\mathbf{M}_{n}\rightarrow I
\]
so that
\[
\psi_{*}(y)=\partial(\varphi_{*}(x)).
\]
The mapping $\varphi\mapsto\psi$ can be chosen to be natural and
well-defined up to homotopy.
\end{thm}

By the ``standard generator'' of $K_{1}\left(\mathbf{I}(\mathbb{K})\right)$
we mean the push-forward of the standard generator of $K_{1}\left(C_{0}(0,1)\right)$
by the obvious inclusion
\[
C_{0}(0,1)
=
C_{0}\left((0,1),\mathbf{M}_{1}\right)\hookrightarrow\mathbf{I}(\mathbb{K}).
\]

\begin{thm}
Let $x$ be the standard generator of $K_{0}\left(\mathbf{M}_{n}\right)$
and $y$ be the standard generator of $K_{1}\left(\mathbf{I}(\mathbb{K})\right).$
Given an ideal $I$ in a $C^{*}$-algebra $A,$ for any $*$-homomorphism\[
\varphi:\mathbb{K}\rightarrow A/I\]
there is a $*$-homomorphism\[
\psi:\mathbf{I}(\mathbb{K})\rightarrow I\]
so that\[
\psi_{*}(y)=\partial(\varphi_{*}(x)).\]

\end{thm}

Notice that from $\psi$ we get realizations of the boundary of $\varphi_{*}(x)$
via a map
\[
\psi_{n}:\mathbf{S}\mathbf{M}_{n}\rightarrow I
\]
for any $n.$

\end{document}